\documentclass[12pt]{article}
\usepackage{amsmath,amssymb,amsbsy}
\usepackage{graphicx}
\usepackage{diagbox,multirow}
\usepackage{tikz}
\usepackage{tkz-graph}
\usetikzlibrary{shapes.geometric}
\usepackage{tkz-euclide,amsmath} 
\usepackage{soul}

\newtheorem{theorem}{Theorem}
\newtheorem{lemma}[theorem]{Lemma}
\newtheorem{corollary}[theorem]{Corollary}
\newtheorem{proposition}[theorem]{Proposition}

\newtheorem{observation}[theorem]{Observation}

\newtheorem{definition}[theorem]{Definition}

\newtheorem{problem}[theorem]{Problem}

\def\gcd{\mathop{\rm gcd}\nolimits}

\newcommand{\qed}{\hfill \rule{.1in}{.1in}}

\def\zet{\mathop\mathbb{Z}\nolimits}
\newcommand{\ZZ}{\mathbb{Z}}

\makeatletter
\def\imod#1{\allowbreak\mkern10mu({\operator@font mod}\,\,#1)}

\title{\bf Group distance magic cubic graphs}%D

\begin{document}

\author{{{Sylwia Cichacz$^{1}$, Štefko Miklavič$^{2,3}$}}
	\footnote{e-mails: cichacz@agh.edu.pl, stefko.miklavic@upr.si}\\
\normalsize $^1$AGH University of Krakow, Poland \\
\normalsize $^2$Andrej Maru\v{s}i\v{c} Institute, University of Primorska, \\
\normalsize Muzejski trg 2, Koper, Slovenia \\
\normalsize $^3$Institute od mathematics, physics and mechanics, \\
\normalsize Jadranska 19, Ljubljana, Slovenia
}

\maketitle
\begin{abstract}
A $\Gamma$\emph{-distance magic labeling} of a graph $G = (V, E)$ with $|V| = n$ is a bijection $\ell$ from $V$ to an Abelian group $\Gamma$ of order $n$, for which there exists $\mu \in \Gamma$, 
such that the weight $w(x) =\sum_{y\in N(x)}\ell(y)$ of every vertex $x \in V$ is equal to $\mu$. In this case, the element $\mu$ is called the \emph{magic constant of} $G$. A graph $G$ is called a \emph{group distance magic} if there exists a $\Gamma$-distance magic labeling of $G$ for every Abelian
group $\Gamma$ of order $n$.

In this paper, we focused on cubic $\Gamma$-distance magic graphs as well as some properties of such graphs.
\end{abstract}
%\begin{keyword}
%Abelian group \sep Constant sum partition \sep group distance magic labeling
%\end{keyword}

%\end{frontmatter}

\section{Introduction}
Throughout this paper $\Gamma$ will denote an Abelian group of order $n$ with the group operation denoted by $+$.  As usually we will write $ka$ to denote $a + a + \cdots + a$ (where the element $a$ appears $k$ times), $-a$ to denote the inverse of $a$ and
we will use $a - b$ instead of $a+(-b)$.  Moreover, the notation $\sum_{a\in S}{a}$ will denote the sum of all elements of the set $S \subseteq \Gamma$.  The identity element of $\Gamma$ will be denoted by $0$. Recall that any group element $\iota\in\Gamma$ of order 2 (i.e., $\iota\neq 0$ such that $2\iota=0$) is called an \emph{involution}.  Let us denote the set of involutions of $\Gamma$ by $I(\Gamma)$.

Let $G=(V,E)$ denote a graph with vertex set $V$ and edge set $E$. For $x \in V$ we denote the set of all neighbors of $x$ by $N(y)$. If there exists an integer $k$ such that $|N(x)|=k$ for every $x \in V$, then we say that $G$ is {\em regular with valency} $k$. Furthermore, we say that $G$ is {\em odd} ({\em even}) {\em regular}, if it is regular with valency $k$, where $k$ is odd (even, respectively). We say that $G$ is {\em cubic} if it is regular with valency $3$.

In \cite{Fro}, Froncek defined the notion of group distance magic graphs, that we now recall. A $\Gamma$\emph{-distance magic labeling} of a graph $G = (V, E)$ with $|V| = n$ is a bijection $\ell$ from $V$ to an Abelian group $\Gamma$ of order $n$, for which there exists $\mu \in \Gamma$, 
such that the weight $w(x) =\sum_{y\in N(x)}\ell(y)$ of every vertex $x \in V$ is equal to $\mu$. In this case, the element $\mu$ is called the \emph{magic constant of} $G$. A graph $G$ is called a \emph{group distance magic}, if there exists a $\Gamma$-distance magic labeling of $G$ for every Abelian
group $\Gamma$ of order $n$.

Concerning cubic graphs, so far only one family of graphs has been studied from the viewpoint of group distance magic labelings, namely $t$ disjoint copies of the complete bipartite graph$K_{3,3}$ \cite{Cic}.
The main result of~\cite{Cic} is even more general, but we will formulate it in this form.

\begin{theorem}[\cite{Cic}]\label{tK_{3,3}} Let $G=tK_{3,3}$ and let $\Gamma$ be an Abelian group of order $6t$. The graph $G$ has a $\Gamma$-distance magic labeling if and only if  $|I(\Gamma)| \ge 2$.
\end{theorem}

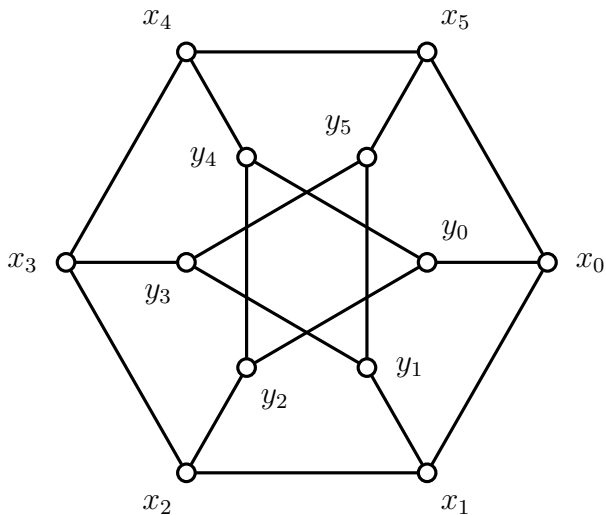
\begin{figure}[ht!]
\begin{center}
\begin{tikzpicture}
  [scale=.8,auto=left,every node/.style={shape=circle,minimum size = 1pt, very thick}]
  \node[draw=black,scale=.6,label=0:$x_0$] (x0) at (4,0) {};
  \node[draw=black,scale=.6,label=300:$x_1$] (x1) at (2,-3.5) {};
  \node[draw=black,scale=.6,label=240:$x_2$] (x2) at (-2,-3.5) {};
  \node[draw=black,scale=.6,label=180:$x_3$] (x3) at (-4,0) {};
  \node[draw=black,scale=.6,label=120:$x_4$] (x4) at (-2,3.5) {};
  \node[draw=black,scale=.6,label=60:$x_5$] (x5) at (2,3.5) {};
  \node[draw=black,scale=.6,label=60:$y_0$] (y0) at (2,0) {};
  \node[draw=black,scale=.6,label=0:$y_1$] (y1) at (1,-1.75) {};
  \node[draw=black,scale=.6,label=300:$y_2$] (y2) at (-1,-1.75) {};
  \node[draw=black,scale=.6,label=240:$y_3$] (y3) at (-2,0) {};
  \node[draw=black,scale=.6,label=180:$y_4$] (y4) at (-1,1.75) {};
  \node[draw=black,scale=.6,label=120:$y_5$] (y5) at (1,1.75) {};
  
  \foreach \from/\to in {x0/x1,x1/x2,x2/x3,x3/x4,x4/x5,x5/x0,
  y0/y2,y1/y3,y2/y4,y3/y5,y4/y0,y5/y1,
  x0/y0,x1/y1,x2/y2,x3/y3,x4/y4,x5/y5}
  \draw[-,very thick] (\from) -- (\to);
  
\end{tikzpicture}
\end{center}
\caption{The generalized Petersen graph $GP(6,2)$.}
\label{generalizedPetersen}
\end{figure}

In this paper, we will consider mostly cubic graphs. In Section~\ref{sec:gp} we consider generalized Petersen graphs. The generalized Petersen graphs were introduced by Coxeter in \cite{ref_Cox} and named by Watkins \cite{ref_Wat}. They are defined as follows. Let $n,k$ be positive integers with $1\leq k<n/2$. Then the generalized Petersen graph $GP(n,k)$ is the graph on $2n$ vertices $\{x_i, y_i : i \in \ZZ_n\}$, where the edge set consists of the polygon $\{x_ix_{i+1} : i \in \ZZ_n\}$, the star polygon $\{y_i y_{i+k} : i \in \ZZ_n\}$, and the spokes $\{x_i y_i : i \in \ZZ_n\}$. In this notation, the ordinary Petersen graph is $G(5,2)$. The generalized Petersen graph $GP(6,2)$ is given in Figure \ref{generalizedPetersen}.

In Section~\ref{sec:union} we study a disjoint unions of group distance magic graphs. We will study $\Gamma$-distance magic labelings of $tG$ (i.e. $t$ disjoint copies of $G$) for a given $\Gamma_0$-distance magic graph $G$, where $\Gamma$ is a group of order $t|\Gamma_0|$.

\section{Preliminary results}

In this section, we recall certain results that we will need later in the paper. Recall that a non-trivial
finite group contains at least one involution if and only if the order of the group is even. The fundamental theorem of finite Abelian groups states that a finite Abelian
group $\Gamma$ of order $n$ can be expressed as the direct product of cyclic subgroups of prime-power order. This implies that
$$
\Gamma\cong\zet_{p_1^{\alpha_1}}\oplus\zet_{p_2^{\alpha_2}}\oplus\cdots\oplus\zet_{p_k^{\alpha_k}}\;\;\; \mathrm{where}\;\;\; n = p_1^{\alpha_1}\cdot p_2^{\alpha_2} \cdots p_k^{\alpha_k}
$$
and $p_i$ for $i \in \{1, 2,\ldots,k\}$ are not necessarily distinct primes. This product is unique up to the order of the factors of the direct product. When $p$ is the number of these cyclic components
whose order is a multiple of $2$, then $\Gamma$ has $2^p-1$ involutions.   By the fundamental theorem of Abelian groups, it is easy to see that any Abelian group $\Gamma$ can be factorized as $\Gamma\cong L\times H$,  with $|L|=2^\eta$ for a non-negative integer $\eta$ and $|H|=\rho$ for an odd positive integer $\rho$. Note that in this case $|I(\Gamma)|=|I(L)|$. Note also that the sum of all elements of a group $\Gamma$ is equal to the sum of all involutions of $\Gamma$. Therefore $\sum_{g\in \Gamma}g= \iota$ if $|I(\Gamma)|=1$
 (where $\iota$ is the unique involution of $\Gamma$) and $\sum_{g\in \Gamma}g= 0$ otherwise (see \cite{CN}, Lemma 8).

 If $H$ is a subgroup of $\Gamma$ then we write $H<\Gamma$. For $g \in \Gamma$ we let $\langle g \rangle$ denote the subgroup of $\Gamma$ generated by $g$.   We will start with some basic results and observations regarding group distance magic labelings of odd regular graphs.
 
 \begin{theorem}[\cite{CicFro}]\label{gr:odd} Let $G$ be an $r$-regular graph on $n$
 	vertices, where $r$ is odd. If $\Gamma$ is an Abelian group of order $n$ with $|I(\Gamma)|=1$, then $G$ is not $\Gamma$-distance magic
 \end{theorem}
 
 \begin{observation}[\cite{CicFro}]\label{obs:odd} Let $G$ be an $r$-regular  graph on $n\equiv 2 \pmod4$ vertices, where $r$ is odd. Then $G$ is not $\Gamma$-distance magic for any Abelian group $\Gamma$ of order $n$.
 \end{observation}

\begin{observation}\label{zero}
Let $G$ be an $r$-regular graph on $n$
vertices, where $r$ is odd. Assume that $\Gamma=A\oplus B$, where $\gcd(|A|,r)=1$ and $|A|\cdot|B|=|\Gamma|=n$. Then the following (i), (ii) are equivalent:
\begin{itemize}
	\item[(i)]
	there exists a $\Gamma$-distance magic labeling of $G$ with a magic constant $(0,\mu_2)\in \Gamma$;
	\item[(ii)]
	there exists a $\Gamma$-distance magic labeling of $\Gamma$ with a magic constant $(\mu_1,\mu_2)\in \Gamma$
	for any $\mu_1\in A$.
\end{itemize} 
\end{observation}
\textit{Proof.} Note that $h\colon A \to A$ defined as $h(g)=rg$, $g\in A$ is an automorphism since $\gcd(|A|,r)=1$.

\noindent
$(i) \Rightarrow (ii)$:
Suppose that $\ell\colon V(G)\to\Gamma$ is a a $\Gamma$-distance magic labeling of $G$ with a magic constant $(0,\mu_2)\in \Gamma$. Pick $\mu_1 \in A$ and let $\ell'\colon V(G)\to\Gamma$
be defined as $\ell'(x)=\ell(x)+(h^{-1}(\mu_1),0)$. Thus $\ell'$ is a bijection such that 
$$
\sum_{y\in N(x)}\ell'(y)=\sum_{y\in N(x)}(\ell(x)+(h^{-1}(\mu_1),0)=(0,\mu_2)+r(h^{-1}(\mu_1),0)=(\mu_1,\mu_2).
$$
Therefore, $\ell'$ is  a $\Gamma$-distance magic labeling of $G$ with a magic constant $(\mu_1,\mu_2)\in \Gamma$.

\noindent
$(ii) \Rightarrow (i)$:
Let $\ell\colon V(G)\to\Gamma$ denote a $\Gamma$-distance magic labeling of $G$ with a magic constant $(\mu_1,\mu_2)\in \Gamma$.  Observe that $\ell'\colon V(G)\to\Gamma$
defined as $\ell'(x)=\ell(x)-(h^{-1}(\mu_1),0)$ is a bijection such that 
$$
\sum_{y\in N(x)}\ell'(y)=\sum_{y\in N(x)}(\ell(x)-(h^{-1}(\mu_1),0)=(\mu_1,\mu_2)-r(h^{-1}(\mu_1),0)=(0,\mu_2).
$$
Therefore, $\ell'$ is  a $\Gamma$-distance magic labeling of $G$ with a magic constant $(0,\mu_2)\in \Gamma$. \qed

 As a corollary, we obtain the following result.
 
\begin{corollary}\label{pof2}
Let $G$ be a cubic graph of order $2^m$ for some $m$. Then there does not exist a $(\zet_2)^m$-distance magic labeling of $G$.
\end{corollary}
\textit{Proof.} Let $\Gamma=(\zet_2)^m$. Suppose that there exists a $\Gamma$-distance magic labeling $\ell$ of $G$.
By Observation~\ref{zero} we can assume $w(x)=0$ for any $x\in V(G)$.
Let $y\in V(G)$ be such that $\ell(y)=0$. Let $z\in N(y)$, and $y_1,y_2\in N(z)$ be the other two neighbors of $z$. Tthen since $w(z)=0=\ell(y)+\ell(y_1)+\ell(y_2)$, we obtain $\ell(y_1)=\ell(y_2)$, a contradiction~\qed

\section{Generalized Petersen graphs}\label{sec:gp}
In this section we will study group distance magic labelings of generalized Petersen graphs $GP(n,k)$ with $k \in \{1,2\}$. We start with some basic observations.

\begin{observation}
Let $G$ be a connected cubic graph on $n$ vertices and let $\Gamma$ be an Abelian group of order $n$. If $G$ is $\Gamma$-distance magic, then there are no $4$-cycles in $G$.\label{4cycle}
\end{observation}
\textit{Proof.} Let $\ell\colon V(G)\to \Gamma$ be a $\Gamma$-distance magic labeling and assume that $v_1v_2v_3v_4v_1$ is a 4-cycle in $G$. If there exists $x\in N(v_1)\setminus N(v_3)$, then there is also $y\in N(v_3)\setminus N(v_1)$ and $\ell(x)=\ell(y)$, a contradiction. This shows that any two antipodal vertices on any 4-cycle of $G$ have equal neighborhoods. In particular, $N(v_1)=N(v_3)=\{v_2,v_4,x\}$  and $N(v_2)=N(v_4)=\{v_1,v_3,z\}$ for some vertices $x,z$ of $G$. Note that $x \ne z$, as otherwise, $x$ would have degree at least $4$, a contradiction. 

Consider now the 4-cycle $v_1v_2v_3xv_1$ of  $G$. By the comment above, $N(v_2) = N(x)$, and so $x$ and $y$ are adjacent. Consequently $G$ is isomorphic to the complete bipartite graph $K_{3,3}$, contradicting Theorem \ref{tK_{3,3}}.\qed

\bigskip

The above observation implies the following result.
\begin{observation}
Let $G=GP(n,1) \;(n \ge 3)$ and let $\Gamma$ denote an Abelian group of order $2n$. Then $G$ is not $\Gamma$-distance magic.
\end{observation}
\textit{Proof.} Note that $G$ contains a $4$-cycle $x_0 y_0 y_1 x_1 x_0$. The result now follows from  Observation~\ref{4cycle}.~\qed

\medskip
To study group distance magic labelings of generalized Petersen graphs $GP(n,2)$, we now introduce a new type of a binary labeling, which we will find useful later.

\begin{definition}
    Let $G$ be a graph with vertex set $V$, $|V|=2n$, and let $\ell: V \to \mathbb{Z}_2$. For each vertex $v \in V$ define the {\em weight of} $v$ by $w(v)=\sum_{u \in N(v)} \ell(u)$. Then $\ell$ is called a {\em balanced zero-neighborhood labeling} of $G$, if $|\ell^{-1}(0)|=|\ell^{-1}(1)|=n$ and $w(v)=0$ for every $v \in V$.
\end{definition}
The following observation is straightforward.

\begin{observation}
\label{bl:ob1}
    Let $G$ be a graph with vertex set $V$, $|V|=2n$, and let $\Gamma$ be a group of order $n$. Assume that $G$ is $\mathbb{Z}_2 \oplus \Gamma$-distance magic, such that the magic constant is of the form $(0,g)$ for some $g \in \Gamma$. Let $\ell: V \to \mathbb{Z}_2 \oplus \Gamma$ be the corresponding $\mathbb{Z}_2 \oplus \Gamma$-distance magic labeling. Let $\ell_1$ ($\ell_2$, respectively) be the projection of $\ell$ on the first (second, respectively) coordinate, that is, for every $v \in V$ we have that $\ell(v)=(\ell_1(v), \ell_2(v))$. Then $\ell_1$ is a balanced zero-neighborhood labeling of $G$.
\end{observation}

\begin{lemma}
\label{bl:lem1}
    Let $\ell$ denote a balanced zero-neighborhood labeling of $GP(n,2)$. Then for every $i \in \mathbb{Z}_n$ we have 
    $$
    \ell(y_i) = \ell(x_{i-1}) + \ell(x_{i+1}), \qquad \ell(x_i) = \ell(y_{i-2}) + \ell(y_{i+2}).
    $$
    In particular, 
    $$
    \ell(x_i) = \ell(x_{i-3}) + \ell(x_{i-1}) + \ell(x_{i+1}) + \ell(x_{i+3}).
    $$
\end{lemma}
\textit{Proof.} Note that $N(x_i)= \{x_{i-1}, x_{i+1}, y_i\}$ and $N(y_i)=\{x_i, y_{i-2}, y_{i+2}\}$. The first part of the lemma now follows immediately from the definition of a balanced zero-neighborhood labeling. The second part of the lemma follows immediately from the first part. \qed

\begin{lemma}
\label{bl:lem2}
    Let $\ell$ denote a balanced zero-neighborhood labeling of $GP(n,2)$. Then for every $i \in \mathbb{Z}_n$ we have 
    $$
    \ell(x_i) + \ell(x_{i+2}) = \ell(x_{i-3}) + \ell(x_{i+5}).
    $$
\end{lemma}
\textit{Proof.} By Lemma \ref{bl:lem1} we have 
$$
\ell(x_i) = \ell(x_{i-3}) + \ell(x_{i-1}) + \ell(x_{i+1}) + \ell(x_{i+3})
$$
and
$$
\ell(x_{i+2}) = \ell(x_{i-1}) + \ell(x_{i+1}) + \ell(x_{i+3}) + \ell(x_{i+5}).
$$
Now add the above two equalities to get the desired result. \qed

\begin{lemma}
\label{bl:lem3}
    Let $\ell$ denote a balanced zero-neighborhood labeling of $GP(n,2)$. Then for every $i \in \mathbb{Z}_n$ and every non-negative integer $t$ we have 
    $$
    \ell(x_i) + \ell(x_{i+3}) = \ell(x_{i+5t}) + \ell(x_{i+5t+3})
    $$
    and 
    $$
    \ell(x_i) + \ell(x_{i+5}) = \ell(x_{i+3t}) + \ell(x_{i+3t+5}).
    $$
\end{lemma}
\textit{Proof.} We prove the first equality. The second one is proved analogously. 
Note that it follows from Lemma \ref{bl:lem2} (replacing index $i$ by $i+3$) that $\ell(x_i) + \ell(x_{i+3}) = \ell(x_{i+5}) + \ell(x_{i+8})$. Replacing index $i$ in this equality by $i+5(t-1)$ we get 
$$
\ell(x_{i+5(t-1)}) + \ell(x_{i+5(t-1)+3}) = \ell(x_{i+5t}) + \ell(x_{i+5t+3}).
$$
The result follows.
\qed

\begin{lemma}
\label{bl:lem4}
    Let $\ell$ denote a balanced zero-neighborhood labeling of $GP(n,2)$. Then for every $i \in \mathbb{Z}_n$ we have 
    $$
    \ell(x_{i+3}) + \ell(x_{i+5}) = \ell(x_{i+18}) + \ell(x_{i+20}).
    $$
\end{lemma}
\textit{Proof.} Using Lemma \ref{bl:lem3} we that

\begin{equation}\label{eq1} 
\begin{split}
\ell(x_i) + \ell(x_{i+3}) = \ell(x_{i+15}) + \ell(x_{i+18}).
\end{split}
\end{equation}

    and 
 \begin{equation}\label{eq2} 
\begin{split}
   \ell(x_i) + \ell(x_{i+5}) = \ell(x_{i+15}) + \ell(x_{i+20}).
\end{split}
\end{equation}   
   
Now add these two equalities to get the desired result.
\qed

\begin{lemma}
\label{bl:lem5}
    Let $\ell$ denote a balanced zero-neighborhood labeling of $GP(n,2)$. Then for every $i \in \mathbb{Z}_n$ we have 
    $$
    \ell(y_i)=\ell(y_{i+15})
    $$
    and 
    $$
    \ell(x_i)=\ell(x_{i+15}).
    $$
\end{lemma}
\textit{Proof.} Combining Lemma \ref{bl:lem1} and Lemma \ref{bl:lem4} we get
$$
\ell(y_i)=\ell(x_{i-1})+\ell(x_{i+1}) = \ell(x_{i+14}) + \ell(x_{i+16}) = \ell(y_{i+15}).
$$
Combining Lemma \ref{bl:lem1} and the above result we also have
$$
\ell(x_i)=\ell(y_{i-2}) + \ell(y_{i+2}) = \ell(y_{i+13}) + \ell(y_{i+17})=\ell(x_{i+15}).
$$
\qed

\begin{definition}
    Let $\ell$ be a balanced zero-neighborhood labeling of a cubic graph $G$ and let $v$ denote a vertex of $G$. We say that $v$ is {\em of type 1} if the labels of all three neighbors of $v$ are $0$. We say that $v$ is {\em of type 2}, if the label of exactly one neighbor of $v$ is $0$. 
\end{definition}

\begin{lemma}
\label{bl:lem6}
    Let $\ell$ denote a balanced zero-neighborhood labeling of $GP(n,2)$. Then for every $i \in \mathbb{Z}_n$ the following (i), (ii) hold:
    \begin{itemize}
        \item[(i)] $x_i$ is of type $2$ if and only if $x_{i+15}$ is of type $2$.
        \item[(ii)] $y_i$ is of type $2$ if and only if $y_{i+15}$ is of type $2$.
    \end{itemize}
    In particular, the number of type $2$ vertices of $GP(n,2)$ is divisible by $n/\gcd(n,15)$.
\end{lemma}
\textit{Proof.} We prove part (i) of Lemma \ref{bl:lem6}. The proof of part (ii) is analogous. Observe that $x_i$ is of type 2 if and only if there exists a neighbor $z$ of $x_i$ with $\ell(z)=1$. By Lemma \ref{bl:lem5} and the definition of generalized Petersen graphs, this will happen if and only if there exists a neighbor $z'$ of $x_{i+15}$ with $\ell(z')=1$, that is, if and only if $x_{i+15}$ is of type 2.

To prove the last part of the above lemma, let $H$ denote the subgroup of $\mathbb{Z}_n$, generated by $\gcd(n,15)$. Note that (i), (ii) above imply that if $x_i$ ($y_i$, respectively) is of type $2$ and $j$ belongs to the same coset of $H$ as $i$, then also $x_j$ ($y_j$, respectively) is of type $2$. As the order of $H$ is $n/\gcd(n,15)$, the result follows.
\qed

\begin{lemma}
\label{bl:lem7}
    Let $\ell$ denote a balanced zero-neighborhood labeling of a cubic graph $G$ on $2n$ vertices. Let $\alpha$ denote the number of type $2$ vertices of $G$. Then $\alpha=3n/2$. In particular, $n$ is even.
\end{lemma}
\textit{Proof.} For the sake of this proof, we treat the labels of vertices of $G$ as integers (and not as elements of $\mathbb{Z}_2$). Let $V$ denote the vertex set of $G$ and note that we have 
$$
\sum_{v \in V} \ell(v) = n.
$$
On the other hand side, the above sum is also equal to
$$
\frac{1}{3} \sum_{v \in V} w(v) = \frac{1}{3} 2 \alpha.
$$
The result follows.
\qed

\begin{proposition}
\label{bl:prop1}
    Let $n \ge 3$. Then $GP(n,2)$ $G$ has no balanced zero-neighborhood labeling.
\end{proposition}
\textit{Proof.} 
Assume to the contrary that $\ell$ is a balanced zero-neighborhood labeling of $G$. Let $\alpha$ denote the number of type $2$ vertices of $G$. Then it follows from Lemma \ref{bl:lem6} and Lemma \ref{bl:lem7} that 
$$
\alpha = \frac{3n}{2} = t \frac{n}{\gcd(n,15)}
$$
for some positive integer $t$. It follows that 
$$
\gcd(n,15) = \frac{2t}{3},
$$
implying that $\gcd(n,15)$ is even, a contradiction.
\qed

\begin{corollary}
	\label{cor:noZ2}
    Let $n \ge 3$. Then $GP(n,2)$ is not $\mathbb{Z}_2 \oplus \Gamma$-distance magic for any group $\Gamma$ or order $n$.
\end{corollary}
\textit{Proof.}  This follows immediately from Observation \ref{bl:ob1} and Proposition \ref{bl:prop1}. \qed

\medskip
Let $\Gamma$ be an Abelian group of order $2n$ and assume that 
$$
\Gamma\cong\zet_{p_1^{\alpha_1}}\oplus\zet_{p_2^{\alpha_2}}\oplus\cdots\oplus\zet_{p_k^{\alpha_k}}\;\;\; \mathrm{where}\;\;\; 2n = p_1^{\alpha_1}\cdot p_2^{\alpha_2} \cdots p_k^{\alpha_k}
$$
The above corollary shows that if $p_i=2$ and $\alpha_i=1$ for some $1 \le i \le k$, then there are no $\Gamma$-distance magic labelings of $GP(n,2)$. When this is not the case, we need a slightly different approach.

\begin{lemma}\label{equations}
	Let $G=GP(n,2)$ and let $\Gamma$ be an Abelian group of order $2n$. If there exists a $\Gamma$-distance magic labeling $\ell$ of $G$, then for any $i\in \ZZ_n$ the following (i), (ii) holds.
	\begin{itemize}
		\item[(i)]
		$$
		\ell(y_{i})+\ell(y_{i+5})=\ell(y_{i+3})+\ell(y_{i+8});
		$$
		\item[(ii)]
		$$
		\ell(x_{i})+\ell(x_{i+5})=\ell(x_{i+3})+\ell(x_{i+8}).
		$$
	\end{itemize}
\end{lemma}
\textit{Proof.} Assume that $\ell$ is a $\Gamma$-distance magic labeling of $G$. We will prove part (i) of the above lemma. A proof for part (ii) is analogous. Since for any $i \in \ZZ_n$ we have $w(y_0)=w(y_i)$, we obtain that 
\begin{equation} \label{eq1}
	\ell(y_{n-2})+\ell(y_2)+\ell(x_0)=\ell(y_{i-2})+\ell(y_{i+2})+\ell(x_i).
\end{equation}
It follows that
\begin{equation} \label{eq2}
	\ell(x_i)=\ell(y_{n-2})+\ell(y_2)-(\ell(y_{i-2})+\ell(y_{i+2}))+\ell(x_0).
\end{equation}
Similarly, using $w(x_0)=w(x_i)$, we obtain that 
\begin{equation} \label{eq3}
	\ell(x_{n-1})+\ell(x_1)+\ell(y_0)=\ell(x_{i-1})+\ell(x_{i+1})+\ell(y_i).
\end{equation}
Applying equation (\ref{eq2}) to equation (\ref{eq3}) we obtain
\begin{equation*} 
	\begin{split}
		\ell(y_{n-2})+\ell(y_2)-(\ell(y_{n-3})+\ell(y_{1}))+\ell(x_0)+\\
		\ell(y_{n-2})+\ell(y_2)-(\ell(y_{n-1})+\ell(y_{3}))+\ell(x_0)+\ell(y_0)=\\
		\ell(y_{n-2})+\ell(y_2)-(\ell(y_{i-3})+\ell(y_{i+1}))+\ell(x_0)+\\
		\ell(y_{n-2})+\ell(y_2)-(\ell(y_{i-1})+\ell(y_{i+3}))+\ell(x_0)+\ell(y_i).
	\end{split}
\end{equation*}
Therefore for any $i\in\ZZ_n$ we have 
\begin{equation*}
	\begin{split}
		\ell(y_{n-3})+\ell(y_{1})+\ell(y_{n-1})+\ell(y_{3})-\ell(y_0)=\\
		\ell(y_{i-3})+\ell(y_{i+1})+\ell(y_{i-1})+\ell(y_{i+3})-\ell(y_i).
	\end{split}
\end{equation*}
Replacing index $i$ in the above equation by $i+5$ and $i+3$ respectively, we get
\begin{equation*} 
	\begin{split}
		\ell(y_{i+2})+\ell(y_{i+6})+\ell(y_{i+4})+\ell(y_{i+8})-\ell(y_{i+5})=\\
		\ell(y_i)+\ell(y_{i+4})+\ell(y_{i+2})+\ell(y_{i+6})-\ell(y_{i+3}),
	\end{split}
\end{equation*}
what implies
\begin{equation}\label{eq4} 
	\begin{split}
		\ell(y_{i})+\ell(y_{i+5})=\ell(y_{i+3})+\ell(y_{i+8}).
	\end{split}
\end{equation}
~\qed

\begin{lemma}
	Let $n \ge 3$. Then $GP(n,2)$ is not $\Gamma$-distance magic for any group $\Gamma$ or order $2n$.
\end{lemma}
\textit{Proof.} Let $\Gamma$ denote an Abelian group of order $2n$ and assume that $\ell$ is a  $\Gamma$-distance magic labeling of $G$. By Lemma~\ref{equations} we have
\begin{equation}\label{eq7} 
	\begin{split}
		\ell(x_{0})+\ell(x_{5}) & =\ell(x_{3j})+\ell(x_{5+3j}),\\
		\ell(x_{1})+\ell(x_{6}) & =\ell(x_{1+3j})+\ell(x_{6+3j}),\\
		\ell(x_{2})+\ell(x_{7}) & =\ell(x_{2+3j})+\ell(x_{7+3j}).
	\end{split}
\end{equation}

For convenience let us denote $\ell(x_0)=a$, $\ell(x_5)=b$, $\ell(x_{10})=c$, and $\ell(x_{15})=d$. By equations~(\ref{eq7}) we have
\begin{equation*}
	\begin{split}
		a+b &=\ell(x_{3j})+\ell(x_{5+3j}),\\
		c+d &=\ell(x_{1+3j})+\ell(x_{6+3j}),\\
		b+c &=\ell(x_{2+3j})+\ell(x_{7+3j}),
	\end{split}
\end{equation*}
for any integer $j$. Thus
\begin{equation*} 
	\begin{split}
		\ell(x_{20}) & =a+b-d,\\
		\ell(x_{25}) & =c+d-a,\\
		\ell(x_{30}) & =a.
	\end{split}
\end{equation*}
Since $\ell$ is bijective, we therefore have  $30 \equiv 0\pmod{n}$, and so $n\in\{3,5,6,10,15,30\}$. By Observation~\ref{obs:odd} we have $n\in\{6,10,30\}$. But then either $\Gamma\cong \ZZ_2\oplus H$ for some group of order $n$, or  $\Gamma\cong \ZZ_4\oplus H$ for some group of odd order $n/2$. However, these two possibilities contradict either Corollary \ref{cor:noZ2} or Theorem \ref{gr:odd}. This finishes the proof.~\qed

\section{Union of group distance magic graphs}\label{sec:union}

In \cite{Wal}, Marr and Wallis defined a Kotzig array as a $p \times k$ grid where each row is a permutation of ${1, 2, \dots, k}$, and each column has a constant sum of $p(k+1)/2$. Notably, a Latin square is a special case of a Kotzig array. Furthermore, a Kotzig array exhibits a "magic-like" property concerning the sums of rows and columns while allowing repeated entries in columns. Kotzig arrays are significant in the context of graph labelings. %There are many constructions based on Kotzig arrays of various magic-type constructions \cite{ref_CicFroSin, ref_CicNik,Fro2,ref_InaLlaMor,ref_McQ,ref_Wal}.

\begin{theorem}[\cite{Wal}]\label{Wal}
A Kotzig array of size $p \times k$ exists if and only if
$p>1$ and $p(k - 1)$ is even. 
\end{theorem}

In \cite{Cic}, the following generalization of Kotzig arrays was introduced. For an Abelian group $\Gamma$ of order $k$, a $\Gamma$-\textit{Kotzig array} of size $p \times k$ is defined as a $p \times k$ grid, where each row is a permutation of the elements of $\Gamma$, and each column has the same sum. Note that if there exists a $\Gamma$-Kotzig array of size $p \times k$, then we could assume without loss of generality that all elements in the first column are equal to $0$ (and so all column sums are equal to $0$). In the same paper, also the following necessary and sufficient conditions for the existence of a $\Gamma$-Kotzig array of size $p \times k$ was provided. 
	
\begin{theorem}[\cite{Cic}]\label{Kotzig}
A $\Gamma$-Kotzig array of size $p \times k$ exists if and only if
$p>1$ and $p$ is even or $|I(\Gamma)|\neq 1$. 
\end{theorem}

\begin{figure}[ht!]
\begin{center}
\begin{tikzpicture}
  [scale=.8,auto=left,every node/.style={shape=circle,minimum size = 1pt, very thick}]
  \node[draw=black,scale=.6] (x0) at (0,4) {};
  \node[draw=black,scale=.6,label=left:{$(1,1,1)$}] (x1) at (-2.57,3.06) {};
  \node[draw=black,scale=.6,label=left:{$(0,1,0)$}] (x2) at (-3.94,0.69) {};
  \node[draw=black,scale=.6,label=left:{$(0,0,0)$}] (x3) at (-3.46,-2.0) {};
  \node[draw=black,scale=.6,label=left:{$(1,1,0)$}] (x4) at (-1.37,-3.76) {};
  \node[draw=black,scale=.6,label=right:{$(0,1,2)$}] (x5) at (1.37,-3.76) {};
  \node[draw=black,scale=.6,label=right:{$(0,0,2)$}] (x6) at (3.46,-2) {};
  \node[draw=black,scale=.6,label=right:{$(1,1,2)$}] (x7) at (3.94,0.69) {};
  \node[draw=black,scale=.6,label=right:{$(0,1,1)$}] (x8) at (2.57,3.06) {};
  \node[draw=black,scale=.6] (y0) at (0,2) {};
  \node[draw=black,scale=.6] (y1) at (-1.73,-1) {};
  \node[draw=black,scale=.6] (y2) at (1.73,-1) {};
  
  \foreach \from/\to in {x0/x1,x1/x2,x2/x3,x3/x4,x4/x5,x5/x6,x6/x7,
  x7/x8,x8/x0,y0/y1,y1/y2,y2/y0,x0/y0,x3/y1,x6/y2,x8/x4,x5/x1,x2/x7}
    \draw[-,very thick] (\from) -- (\to);
  \draw (0,4.5) node {$(0,0,1)$};
\draw (1,2.5) node {$(1,0,1)$};
\draw (2.75,-0.75) node {$(1,0,2)$};
\draw (-2.75,-0.75) node {$(1,0,0)$};

\end{tikzpicture}
\end{center}
\caption{$\zet_2\oplus\zet_2\oplus\zet_3$-distance magic labeling of Tietze graph.}
\label{Tietze}
\end{figure}
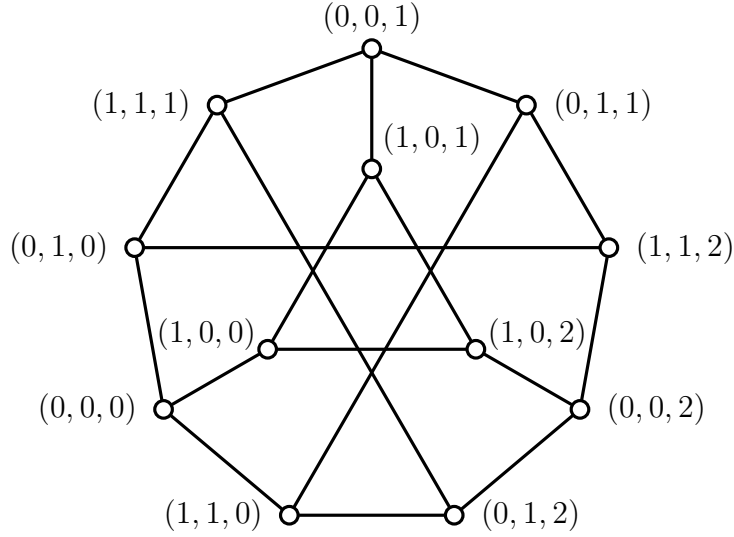

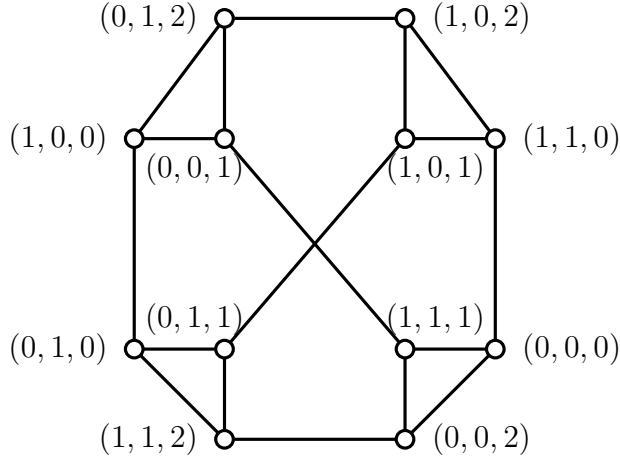
\begin{figure}[ht!]
\begin{center}
\begin{tikzpicture}
  [scale=.8,auto=left,every node/.style={shape=circle,minimum size = 1pt, very thick}]
  \node[draw=black,scale=.6,label=right:{$(1,0,2)$}] (x0) at (1.5,3) {};
  \node[draw=black,scale=.6,label=left:{$(0,1,2)$}] (x1) at (-1.5,3) {};
  \node[draw=black,scale=.6,label=left:{$(1,0,0)$}] (x2) at (-3,1) {};
  \node[draw=black,scale=.6,label=left:{$(0,1,0)$}] (x3) at (-3,-2.5) {};
  \node[draw=black,scale=.6,label=left:{$(1,1,2)$}] (x4) at (-1.5,-4) {};
    \node[draw=black,scale=.6,label=right:{$(0,0,2)$}] (x5) at (1.5,-4) {};
  \node[draw=black,scale=.6,label=right:{$(0,0,0)$}] (x6) at (3,-2.5) {};
  \node[draw=black,scale=.6,label=right:{$(1,1,0)$}] (x7) at (3,1) {};

  \node[draw=black,scale=.6] (y0) at (-1.5,1
) {};

  \node[draw=black,scale=.6] (y1) at (-1.5,-2.5) {};

    \node[draw=black,scale=.6] (y2) at (1.5,-2.5) {};
 
\node[draw=black,scale=.6] (y3) at (1.5,1) {};
 % \node[draw=black,scale=.6] (y0) at (0,2) {};
 % \node[draw=black,scale=.6] (y1) at (-1.73,-1) {};
  %\node[draw=black,scale=.6] (y2) at (1.73,-1) {};
  
  \foreach \from/\to in {x0/x1,x1/x2,x2/x3,x3/x4,x4/x5,x5/x6,x6/x7,
  x7/x0,x1/y0,x2/y0,y0/y2,x3/y1,x4/y1,y1/y3,x5/y2,x6/y2,x7/y3,x0/y3}
    \draw[-,very thick] (\from) -- (\to);
  \draw (-2,0.5) node {$(0,0,1)$};
\draw (2,0.5) node {$(1,0,1)$};
\draw (-2,-2) node {$(0,1,1)$};
\draw (2,-2) node {$(1,1,1)$};

\end{tikzpicture}
\end{center}
\caption{$\zet_3\oplus\zet_2\oplus\zet_2$-distance magic labeling of a graph $G$ }
\label{ppgraph}
\end{figure}

Let $G=(V,E)$ denote a graph with vertex-set $V$. For a positive integer $t$, we denote by $tG$ the disjoint union of $t$ copies of $G$. Let  $V_1, V_2, \ldots ,  V_p$ be a partition of $V$  (that is, $V = V_1\cup V_2\cup\ldots\cup V_p$ where $V_i\cap V_j = \emptyset$ for $i\neq j$ and $V_i \ne \emptyset$ for $1 \le i \le p$). Following \cite{CicNik}, we say that $V_1, V_2, \ldots ,  V_p$ form a $p$-\emph{partition} of $G$, if for every $x \in V(G)$ we have
\[ |N(x)\cap V_1|=|N(x)\cap V_2|=\cdots=|N(x)\cap V_p|.\]
Note that every graph has a $1$-partition (with $V_1=V$), and so we assume that $p \ge 2$. Observe that if $G$ is $r$-regular and has a $p$-partition, then $p$ divides $r$. Thus one can easily see that Tietze graph given in Figure~\ref{Tietze} does not have a $p$-partition for any $p \ge 2$. On the other hand, consider a graph $G$ presented in Figure~\ref{ppgraph}. This graph has a $\zet_2\oplus\zet_2\oplus\zet_3$-distance magic labeling $\ell$, which is also presented in Figure~\ref{ppgraph}. Let $V$ denote the vertex-set of $G$. For $1 \le i \le 3$, let $V_i$ consist of vertices $v \in V$, for which the third coordinate of $\ell(v)$ equals $i-1 \in \zet_3$. It is easy to see that $V_1, V_2, V_3$ form a $3$-partition of $G$.

We now prove two lemmas which will play an important role in the proof of the main result of this section.

\begin{lemma}\label{lemgl}
Let $A,B$ be Abelian groups and let $|B|=t$.  Let $G$ be an $r$-regular graph having a  $p$-partition. If  $G$ is $A$-distance magic and $p$ is even or $|I(B)|\neq 1$, then the graph $tG$ is $\Gamma$-distance magic, where $\Gamma=A \oplus B$.
\end{lemma}
\textit{Proof.} 
Let $V$ denote a vertex-set of $G$ and let  $V_1, V_2, \ldots, V_p$ be a $p$-partition of $G$. Let $G^1,G^2,\ldots, G^t$ denote disjoint copies of $G$, that form the graph $tG$. For a vertex $x\in V$, let $x^j$ denote the corresponding vertex of $G^j$. For $1 \le i \le p$ and $1 \le j \le t$ we define $V_i^j = \{x^j : x \in V_i\}$. Observe that  $V_1^j, V_2^j, \ldots, V_p^j$ is a $p$-partition of $G^j$. 

Let $\ell'$ be an $A$-distance magic labeling of $G$ with the magic constant $\mu'$. 
Let $K_{B}=(k_{i,j})$ be a $B$-Kotzig array  of size $p \times t$ which exists by Theorem~\ref{Kotzig}. Without loss of generality, we assume that the column sums of $K_B$ are equal to $0$. 
Define $\ell \colon V(tG) \to A \oplus B$ as follows. For $x \in V_i \; (1 \le i \le p)$ and $1 \le j \le t$ we let 
$$
\ell(x^j) =(\ell'(x), k_{i,j}). 
$$
Let us first verify that $\ell$ is a bijection. Suppose that for $x \in V_i$ and $y \in V_{i'}$ we have that $\ell(x^j)= \ell(y^{j'})$ for some $1 \le j,j'\le t$. 
Then $(\ell'(x), k_{i,j})=(\ell'(y), k_{i',j'})$, what implies $\ell'(x)=\ell'(y)$. As $\ell'$ is a bijection, we have that $x=y$, and consequently $i=i'$. From this it follows that $k_{i,j}=k_{i,j'}$, and so $j=j'$. This shows that $\ell$ is an injection. Since $|V(tG)|=|A \oplus B|$ is finite, it follows that $\ell$ is also a bijection.

It is now easy to see that $\ell$ is a $\Gamma$-distance magic labeling of $tG$ with a magic constant $\mu = (\mu',0)$. Indeed, pick $x \in V$ and $1 \le j \le t$, and consider the weight of $x^j$. We have 
$$
  w(x^j) = \sum_{y \in N(x)} \ell(y^j) = \sum_{i=1}^p \sum_{y \in N(x) \cap V_i} (\ell'(y),k_{i,j}) = (\mu',\sum_{i=1}^p \sum_{y \in N(x) \cap V_i} k_{i,j}) =
$$
$$
  (\mu',\frac{r}{p} \sum_{i=1}^p k_{i,j}) = (\mu',0). 
$$
This concludes the proof. ~\qed

\medskip
The proof of the lemma below is similar to the proof of Lemma~\ref{lemgl}, but instead of Kotzig array on an Abelian group we will use a regular Kotzig array on integers.

\begin{lemma}\label{lemgl2}
Let $A$ be an Abelian group and let $\zet_{tm}$ be a cyclic group of order $tm \; (t,m \ge 1)$. Let $\Gamma = A \oplus \zet_{tm}$ and $\Gamma_0 = A \oplus \left\langle  t\right\rangle \le \Gamma$. Let $G$ be an $r$-regular graph having a  $p$-partition. If  $G$ is $\Gamma_0$-distance magic and $p(t-1)$ is even, then $tG$ is $\Gamma$-distance magic.
\end{lemma}
\textit{Proof.} 
Let $V$ denote a vertex-set of $G$ and let  $V_1, V_2, \ldots, V_p$ be a $p$-partition of $G$. Let $G^1,G^2,\ldots, G^t$ denote disjoint copies of $G$, that form the graph $tG$. For a vertex $x\in V$, let $x^j$ denote the corresponding vertex of $G^j$. For $1 \le i \le p$ and $1 \le j \le t$ we define $V_i^j = \{x^j : x \in V_i\}$. Observe that  $V_1^j, V_2^j, \ldots, V_p^j$ is a $p$-partition of $G^j$. 

By Theorem~\ref{Wal} there exists a Kotzig array $K=(k_{i,j})$ of size $p\times t$.  Recall that for every $1 \le j \le t$, the sum of the elements of $j$-th column of $K$ is equal to $\sum_{i=1}^p k_{i,j}=p(t+1)/2$. Let $\ell'$ be a $\Gamma_0$-distance magic  labeling of $G$ with the magic constant $\mu'$. For any $v\in V$ let $\ell'(v)=(\ell'_1(v),\ell'_2(v))\in A\oplus\left\langle  t\right\rangle$. 

Define $\ell \colon V(tG) \to A \oplus \zet_{tm}$ as follows. For $x \in V_i \; (1 \le i \le p)$ and $1 \le j \le t$ we let 
$$
  \ell(x^j) =(\ell'_1(x),\ell'_2(x)+ k_{i,j}),
$$
where the sum on the second coordinate is modulo $tm$. Let us first verify that $\ell$ is a bijection. Suppose that for $x \in V_i$ and $y \in V_{i'}$ we have that $\ell(x^j)= \ell(y^{j'})$ for some $1 \le j,j'\le t$. Then $\ell'_1(x)=\ell'_1(y)$ and 
\begin{equation}
	\label{eq:tm}
	\ell'_2(x)+ k_{i,j}=\ell'_2(y)+ k_{i',j'}. 
\end{equation}
Now \eqref{eq:tm} implies that
$$
  k_{i,j}-k_{i',j'} = \ell'_2(y)-\ell'_2(x) \in \langle t \rangle.
$$
This shows that $k_{i,j}-k_{i',j'}$ is divisible by $t$. Since $1 \le k_{i,j}, k_{i',j'} \le t$, this implies $k_{i,j}=k_{i',j'}$, and so \eqref{eq:tm} yields $\ell'_2(x)=\ell'_2(y)$. Therefore $\ell'(x)=\ell'(y)$, and as $\ell'$ is a bijection, we have that $x=y$. This shows that $\ell$ is an injection. Since $|V(tG)|=|A \oplus \zet_{tm}|$ is finite, it follows that $\ell$ is also a bijection.

It remains to show that there exists a magic constant for the labeling $\ell$. This could be done similarly as in the proof of Lemma \ref{lemgl}. We get that $\ell$ is a $\Gamma$-distance magic  labeling of $tG$ with a magic constant 
$$
\mu = \mu'+\Big(0,\frac{r}{p} \cdot \frac{p(t+1)}{2} \Big) = \mu'+\Big( 0, \frac{r(t+1)}{2} \Big).
$$
~\qed

\medskip
We are now ready to state our main result regarding group distance magic labelings of disjoint union of graphs.
\begin{theorem}
	Let $G$ be an $r$-regular graph having a  $p$-partition. Assume that $G$ is $H'$-distance magic for an Abelian group $H'$. Let $\Gamma$ be an Abelian group which contains a subgroup $H \le \Gamma$, such that $H' \cong H$. Let $t=|\Gamma|/|H|$ denote the index of $H$ in $\Gamma$.	Assume that either  $p$ is even or that $t$ is odd. Then $tG$ has a $\Gamma$-distance magic labeling.
\end{theorem}
\textit{Proof.} 
By the Fundamental Theorem of Finite Abelian Groups, the groups $\Gamma$ and $H$ can be written as a direct products
$$
\Gamma\cong\zet_{q_1^{\alpha_1}}\oplus\zet_{q_2^{\alpha_2}}\oplus\cdots\oplus\zet_{q_w^{\alpha_w}},
$$
$$
H \cong \zet_{q_1^{\beta_1}}\oplus\zet_{q_2^{\beta_2}}\oplus\ldots\oplus\zet_{q_w^{\beta_w}},
$$ 
where $q_1,\ldots,q_w$ are (not necessary distinct) primes  and  $0\leq \beta_i\leq \alpha_i \; (1 \le i \le w)$. Note that $t=|\Gamma|/|H|$ implies $t=q_1^{\alpha_1-\beta_1}\cdot q_2^{\alpha_2-\beta_2} \cdots q_w^{\alpha_w-\beta_w}$. For $0 \le i \le w$ define  
$$
\Gamma_i \cong\zet_{q_1^{\alpha_1}} \oplus \zet_{q_2^{\alpha_2}}\oplus\cdots\oplus \zet_{q_i^{\alpha_i}} \oplus \zet_{q_{i+1}^{\beta_{i+1}}} \oplus\cdots\oplus \zet_{q_w^{\beta_w}},
$$
and note that $\Gamma_0 \cong H$ and $\Gamma_w \cong \Gamma$.  For $0 \le i \le w$ we also define
$$
  t_i =  q_1^{\alpha_1-\beta_1}\cdot q_2^{\alpha_2-\beta_2} \cdots q_i^{\alpha_i-\beta_i}.
$$
Note that $t_0=1$ and $t_w=t$. We claim that $t_iG$ has a $\Gamma_i$-distance magic labeling. We will prove our claim by induction on $i$. Observe that the claim is true for $i=0$ since $H' \cong H \cong \Gamma_0$ and $t_0=1$. Assume now that the claim is true for some $0 \le i \le w-1$. We will show that the claim is true also for $i+1$. Assume first that $\beta_{i+1}=0$. Note that if $p$ is odd, then $t$ must be odd by our assumption, and so $q_{i+1}$ is odd. It follows that for $A=\Gamma_i$ and $B=\zet_{q_{i+1}^{\beta_{i+1}}}$ the assumptions of Lemma \ref{lemgl} are satisfied, and so $q_{i+1}^{\beta_{i+1}}(t_iG)=t_{i+1}G$ has a $\Gamma_{i+1}$-distance magic labeling  (note that $\Gamma_{i+1} \cong A \oplus B$).

Assume next that $\beta_{i+1} \ge 1$. If $\beta_{i+1}=\alpha_{i+1}$, then $\Gamma_{i+1}=\Gamma_i$ and $t_{i+1}=t_i$, and so  $t_{i+1}G$ has a $\Gamma_{i+1}$-distance magic labeling by the induction hypothesis. Therefore, assume $\beta_{i+1} < \alpha_{i+1}$. If $p$ is odd, then again our assumption guarantees that $t$ is also odd, and so $q_{i+1}$ is odd too. Now let 
$$
A=\zet_{q_1^{\alpha_1}} \oplus \zet_{q_2^{\alpha_2}}\oplus\cdots\oplus \zet_{q_i^{\alpha_i}} \oplus \zet_{q_{i+2}^{\beta_{i+2}}} \oplus\cdots\oplus \zet_{q_w^{\beta_w}}.
$$
By Lemma \ref{lemgl2} and induction hypothesis, $q_{i+1}^{\alpha_{i+1}-\beta_{i+1}}(t_i G) = t_{i+1} G$ has a $\Gamma_{i+1}$-distance magic labeling.
~\qed\\

Let $t$ be an odd integer. Note that  for  a group $\Gamma$ of order $12t$ such that $|I(\Gamma)|\neq 1$, there is $\zet_2 \oplus \zet_2 \oplus \zet_3\leq\Gamma$. Therefore, the above theorem implies that for the graph $G$ given in Figure~\ref{ppgraph} and $t$ odd, the union $tG$ is $\Gamma$-distance magic for any $\Gamma$ of order $12t$ such that $|I(\Gamma)|\neq1$.

In some cases, $tG$ is $\Gamma$-distance magic even if $G$ does not have a $p$-partition. This follows from our next result.
\begin{lemma}\label{obsgl}
Let $A,B$ be Abelian groups and let us denote $|B|=r$. If  an $r$-regular graph $G$  is $A$-distance magic, then $rG$ is $\Gamma$-distance magic, where $\Gamma = A \oplus B$.
\end{lemma}
\textit{Proof.} 
Let $V$ denote the vertex-set of $G$ and let $G^1,G^2,\ldots, G^r$ denote disjoint copies of $G$, that form the graph $rG$. For a vertex $x\in V$, let $x^j$ denote the corresponding vertex of $G^j$. Let us fix an ordering $b_1, b_2, \ldots, b_r$ of elements of $B$. Let $\ell'$ be an $A$-distance magic  labeling of $G$ with the magic constant $\mu'$. For any $x \in V$ and $1 \le j \le r$ define
$$
\ell(x^j) =(\ell'(x),b_j).
$$
It is easily seen that $\ell$ is a $\Gamma$-distance magic  labeling of the graph $rG$ with a magic constant $\mu = (\mu',0)$.~\qed

The above lemma implies that for the Tietze graph $G$ given in Figure~\ref{Tietze}  the  union $3^mG$ is $(\zet_2)^2\oplus(\zet_3)^{m+1}$-distance magic for any $m\geq 0$.

\section{Conclusions}
In this paper, on one hand we proved that there does not exist a generalized Peteresen graph $GP(n,j)$ for $j\in\{1,2\}$ that admits a $\Gamma$-distance magic labeling, where $\Gamma$ is an Abelian group of order $2n$. On the other hand, we showed a new infinite family of cubic $\Gamma$-distance magic graphs. We also defined a new type of binary labeling called balanced zero-neighborhood labeling. We finish this section with the following open problem.
\begin{problem}
Characterize cubic graphs that admit balanced zero-neighborhood labeling.
\end{problem}

\section{Statements and Declarations}
The work of the first author was partially  supported by program ''Excellence initiative – research university'' for the AGH University. The work of the second author is supported in part by the Slovenian Research and Innovation Agency (research program P1-0285 and research projects J1-3001, J1-3003, J1-4008, J1-4084, N1-0208, N1-0353, J1-50000, J1-60012).


\begin{thebibliography}{10}

%\bibitem{ref_AnhCicPetTep2}  M. Anholcer, S. Cichacz, I. Peterin,    A. Tepeh, \textit{Distance magic labeling and two product graphs}, Graphs and Combinatorics 31(5) (2015), 1125--1136, DOI: 10.1007/s00373-014-1455-8.



%\bibitem{AFK} S. Arumugam, D. Froncek, N. Kamatchi, \emph{Distance Magic Graphs---A Survey}, Journal of the Indonesian Mathematical Society, Special Edition (2011) 11--26.

%\bibitem{C} S. Cichacz, \emph{Distance magic $(r,t)$-hypercycles},  Utilitas Mathematica (2013), accepted.

\bibitem{Biggs}
N. L. Biggs,  \textit{Algebraic Graph Theory}, 2nd ed. Cambridge, England: Cambridge University Press, 1993.

\bibitem{Cic} S. Cichacz, \emph{Zero sum partition of Abelian groups into sets of the same order and its applications}, Electronic Journal of Combinatorics 25(1) (2018), \#P1.20. 

%\bibitem{Cic1} S. Cichacz, \emph{Note on group distance magic graphs $G[C_4]$}, Graphs and Combinatorics 30(3) (2014) 565--571, DOI: 10.1007/s00373-013-1294-z.

%\bibitem{Cic3} S. Cichacz, \emph{On zero sum-partition of Abelian groups into three sets and group distance magic labeling}, Ars Mathematica Contemporanea 13(2) (2017), 417--425. 


\bibitem{CicFro} S. Cichacz, D. Froncek, \emph{Distance magic circulant graphs},  Discrete Mathematics 339(1) (2016) 84-94, DOI: 10.1016/j.disc.2015.07.002.

\bibitem{CicNik}
S. Cichacz, M. Nikodem, \textit{Note on union of distance magic graphs}, Discussiones Mathematicae - Graph Theory 37 (2017) 239-249

\bibitem{CN}D. Combe, A.M. Nelson, W.D. Palmer, \emph{Magic labellings of graphs over finite abelian groups}, Australasian Journal of Combinatorics 29 (2004) 259--271.
\bibitem{ref_Cox}
H. S. M. Coxeter, Self-dual configurations and regular graphs, Bull. Amer. Math. Soc. 56 (1950), 413--455.

\bibitem{Fro} D. Froncek,
\emph{Group distance magic labeling of of Cartesian products of cycles}, Australasian
Journal of Combinatorics 55 (2013) 167--174.

%\bibitem{HR} N. Hartsfield, G. Ringel, Pearls in Graph Theory, pp. 108--109, Academic Press, Boston (1990) (Revised version 1994).

\bibitem{Wal} A.M. Marr, W.D. Wallis, \textit{Magic Graphs}, 2nd ed., Springer, 2013.

%\bibitem{MRS}M. Miller, C. Rodger and R. Simanjuntak, \emph{Distance magic labelings of graphs}, Australasian Journal of Combinatorics, 28 (2003), 305--315.

%\bibitem{ONeSla}A. O'Neal, P.J. Slater, \textit{Uniqueness of vertex magic constants}, SIAM J. Discrete Math. 27 (2013), 708--716.

\bibitem{ref_Wat}
M. E. Watkins (1969), A theorem on tait colorings with an application to the generalized Petersen graphs, Journal of Combinatorial Theory, Volume 6, Issue 2, 152--164, https://doi.org/10.1016/S0021-9800(69)80116-X.


%bibitem{StrSch}H.J. Straight, P. Schillo, \textit{On the problem of partitioning $\{1,2,\ldots,\}$ into subsets having equal sums}, Proc. AMS, 74 (1979),  229--231.

%\bibitem{Rao}S.B. Rao, \emph{Sigma Graphs---A Survey}, In Labelings of Discrete Structures and Applications, eds. B.D. Acharya, S. Arumugam and A. Rosa, Narosa Publishing House, New Delhi, (2008), 135--140.

%\bibitem{Vi}V. Vilfred, \emph{$\Sigma$-labelled Graphs and Circulant Graphs}, Ph.D. Thesis, University of Kerala, Trivandrum, India, 1994.
\end{thebibliography}
\end{document}